\documentclass[reqno,a4paper,11pt]{amsart}
\usepackage{amsmath, amsthm, amscd, amsfonts, amssymb, graphicx, color}
\usepackage{amssymb}
\usepackage[english]{babel}
\usepackage[dvips, lmargin=4cm, rmargin=4cm, tmargin=4cm, bmargin=4cm]{geometry}
\usepackage[colorlinks=true,urlcolor=blue,linkcolor=blue]{hyperref}
\usepackage{hyperref}
\usepackage[T1]{fontenc}
\usepackage[english]{babel}

\usepackage{enumitem}
\usepackage{dsfont}
\newtheorem{thm}{Theorem}[section]
\newtheorem{defini}{Definition}[section]
\newtheorem{rem}{Remark}[section]
\newtheorem{lem}{Lemma}[section]
\newtheorem{prop}{Proposition}[section]
\newtheorem{coro}{Corollary}[section]

\def \N {\mathbb{N} }

\def \R {\mathbb{R} }

\begin{document}
\title[Eigenvalue problems]{Eigenvalue problems in fractional anisotropic Orlicz-Sobolev spaces}
\author[Mohammed SRATI]
{Mohammed SRATI}
\address{ Mohammed SRATI\newline
 High School of Education and Formation (ESEF), University Mohammed First, Oujda,
 Morocco.}
\email{srati93@gmail.com}

\maketitle
\begin{abstract}
In this paper, we introduce the fractional anisotropic Orlicz-Sobolev spaces, and by using some variational methods, we establish the existence or non-existence of eigenvalues of fractional anisotropic problems involving a nonlocal
integro-differential operator of elliptic type.  In each case, the competition between the growth rates of the anisotropic coefficients plays an essential role in the description of the set of eigenvalues. 
\end{abstract}
\section{Introduction}
When we introduce the notion of derivative, we quickly realize that we can apply the concept of derivative to the derivative function itself, and by the same concept we can introduce the second derivative, then the successive derivatives of integer order. Integration, is an inverse operation of the derivative, can possibly be considered as a derivative of the order minus one. We can also wonder if these derivatives of successive orders have a fractional-order equivalent.
Fractional derivation theory is a subject almost as old as classical calculus as we know it today, its origins date back to the late 17th century, the time when Newton and Leibniz developed the foundations of calculus differential and integral.\\

The idea of generalizing the notion of differentiation $\frac{d^nf(x)}{dx^n}$ to noninteger
orders of $n$ appeared at the birth of the differential calculus itself. The
first attempt to discuss such an idea recorded in history was contained in
the correspondence of Leibniz. In one of his letters to Leibniz concerning
the theorem on the differentiation of a product of functions, Bernoulli
asked about the meaning of this theorem in the case of noninteger order
of differentiation. Leibniz in his letters to L'H\"{o}pital (1695) made some
remarks on the possibility of considering differentials and derivatives of
order $1/2$. In 1820 Lacroix showed an exact formula for the evaluation of
the derivative $\frac{d^{1/2}f(x)}{dx^{1/2}}$.
Therefore, the fractional derivation was developed intensely, therefore, different examples of such operators are present in literature. Among these ones, Riemann-Liouville and Caputo fractional derivatives are the most exploited in one-dimensional applications.\\

 After the construction of these different types of fractional derivatives, It was quite natural to ask for spaces $W^{s,p}$ with $s$ fractional, filling the gaps between $L^p$, $W^{1,p}$, $W^{2,p}$,... Several proposals were made in the fifties and these attempts culminated around 1960.  Let us first consider defining the H\"{o}lder spaces $C^s$ with $0 < s\neq$  integer (filling the gaps between the
 spaces $C^k$ with $k\in \N^*$). Let $0 < \rho < 1$, then we introduce the norm
$$||f||_{C^\rho}:=||f||_{C}+\sup\dfrac{|f(x)-f(y)|}{|x-y|^\rho}$$
where the supremum is taken over all $x\in \R^N$ and  $y\in \R^N$  with $x \neq y$. Let $s\in\R$, then we put
  $$s=[s]+\left\lbrace s\right\rbrace $$
  where $[s]$ it is the integer part of $s$ and $0 <\left\lbrace s \right\rbrace <1 $. Then for $ 0 <s \neq$ integer, the H\"{o}lder space is defined as
$$C^s=\left\lbrace f\in C: ||f||_{C^s}=||f||_{C^{[s]}}+\sum_{|\alpha=[s]}||D^\alpha f||_{C^{\left\lbrace s\right\rbrace }}<\infty\right\rbrace. $$

H\"{o}lder spaces have been employed by J.
Schauder and other mathematicians since the mid-thirties in connection with boundary value problems for second order elliptic differential equations. The final step in
this direction is due to $C$. Miranda, see his book \cite{mira} where one finds also many
references. S. Agmon, A. Douglis and L. Nirenberg extended these considerations to
boundary value problems for higher order elliptic differential equations, see \cite{ADN}.\\

The spaces introduced by S.M. Nikol'skij in 1951, see \cite{nik}, can be obtained from the H\"{o}lder space   if one replaces there the $C$-norm by an $L^p$ norm: Let $0 < s \neq$ integer and $1 < p < \infty$, then $B^s_{p,\infty}$ is the collection of all $f\in L^p$ such that
\begin{equation}\label{intro1eg}
  ||f||_{B^s_{p,\infty}}=||f||_{W^{[s]}_p}+\sum_{|\alpha=[s]}\sup\limits_{\R^N\ni h\neq 0}|h|^{-\left\lbrace s\right\rbrace }||\Delta_hD^\alpha f||_{L^p}\end{equation}
  is finite.\\

  N. Aronszajn \cite{aro}, L.N. Slobodeckij \cite{slo} and E. Gagliardo \cite{gag} suggested in 1955-1958 to replace the sup-norm in $(\ref{intro1eg})$ with respect to $h$ by an $L^p$-norm: Let $0 < s\neq$
  integer and $1 < p < \infty$ then $B^s_{p,p}$ s is the collection of all $f\in L^p$ such that
$$||f||_{B^s_{p,p}}=||f||_{W^{[s]}_p}+\sum_{|\alpha=[s]}\left( \int_{\R^N}|h|^{-\left\lbrace s\right\rbrace p  }||\Delta_hD^\alpha f||^p_{L^p}\dfrac{dh}{|h|^N}\right)^\frac{1}{p}.$$
is finite. Let, for example,
$0 < s < 1$, then we get the fractional Sobolev space $ W^{s,p}$ defined as the collection of all $f\in L^p$ such that
 $$||f||_{B^s_{p,p}}=||f||_{L^p}+\left( \int_{\R^N}\int_{\R^N}\dfrac{|f(x)-f(y)|^p}{|x-y|^{sp+N}}dxdy\right)^\frac{1}{p},$$
which makes clear that we replaced the sup-norm in the H\"{o}lder spaces, 
by an $L^p$-norm.\\

The fractional Sobolev spaces have been a classical topic in functional and harmonic analysis all along, and some important books, such as \cite{lan} treat the topic in detail.  Recently, great attention has been focused on the study of fractional spaces, and nonlocal operators of elliptic type, both for pure mathematical research and in view of concrete real-world applications. This type of operator arises in a quite natural way in many different contexts, such as, among others, the thin obstacle problem, optimization, finance, phase transitions, stratified materials, anomalous diffusion, crystal dislocation, soft thin films, semipermeable membranes, flame propagation,
conservation laws, ultrarelativistic limits of quantum mechanics, quasi-geostrophic
flows, multiple scattering, minimal surfaces, materials science, water waves,
chemical reactions of liquids, population dynamics, geophysical fluid dynamics, and
mathematical finance (American options). \\

Even with this considerable interest in (fractional and classical) Lebesgue spaces, but in certain equations, precisely in the case of nonhomogeneous operators when we try to consider certain conditions on these operators, the problem cannot be formulated with the classical  Lebesgue and Sobolev spaces.  As an example of application in the field of image processing, let us consider the restoration model with the following additive noise:
 \begin{equation}\label{int8eg}
  f=u+bruit.
  \end{equation} 
The energy function associated with the equation $ (\ref{int8eg})$ is then:
  \begin{equation}\label{int9eg}
    J(u):=\int_\Omega(f-u)^2dx,
    \end{equation}
    with $\Omega$ a bounded open of $ \R^2 $. As this is a poorly conditioned problem, then we need to add an energy regularization term, which is usually written in the form:
  $$\int_\Omega A(|\nabla u|)dx,$$
  this term represents diffusion. The energy function then takes the following form:
  $$J(u):=\dfrac{\lambda}{2}(f-u)^2dx+\int_\Omega A(|\nabla u|)dx.$$
  The problems of minimization with sublinear regularization terms of this kind are in general badly posed in spaces with bounded variation (BV). On the other hand, thanks to the Orlicz-Sobolev spaces, this model of restoration, not only is well-posed but also the existence and the uniqueness of a minimum are guaranteed (see   \cite[Theorem 3.1.1]{samar}). Hence the interest of the Orlicz spaces. This functional framework is a generalization of the classical Lebesgue spaces $L^p$. They were first introduced by Birnbaum-Orlicz in \cite{birnh} and Orlicz in \cite{orlicz}.\\ 

In this sense, In the last years, great attention has been devoted to introduce the fractional version of the Orlicz-Sobolev spaces, namely a fractional Sobolev space constructed from an Orlicz space instead of $L^p$, who will be named: the fractional Orlicz-Sobolev space. The first proposals to define this space are given by Salort et al \cite{sal1} in 2017, and  some proprieties of these spaces were later developed by some authors : Azroul et al \cite{sr5,3,sr_mo}, Bahrouni et al \cite{bh1,bh2} and Alberico et al \cite{ang1,ang2}. Consequently, the theory of nonlocal (or integro differential) operators in fractional Orlicz-Sobolev  spaces has seen an important development, and several problems involving the  fractional $ a(.)$-Laplacian operators have been studied,  in which the authors have used different methods to get the existence of solutions (see \cite{srati3,SR,SRT,sal2,sal3}).\\


In this paper we propose  to study the following eigenvalue problem 
$$ \label{P}
 (P_a) \hspace*{0.1cm} \left\{ 
   \begin{array}{clclc}
\sum\limits_{i=1}^{N}(-\Delta)^{s}_{a_i(.)} u & = & \lambda |u|^{q(x)-2}u   & \text{ in }& \Omega \\\\
   \hspace*{1cm} u & = & 0 \hspace*{0.2cm} \hspace*{0.2cm} & \text{ in } & \R^N\setminus \Omega,
   \end{array}
   \right. 
$$
 where $\Omega$ is a Lipschitz open bounded subset of $\mathbb{R}^N$, $N\geqslant 1$, $i\in\{1,...,N\}$ $q:\overline{\Omega}\rightarrow(1,+\infty)$ is a  bounded continuous function, $0<s<1$ and $\lambda$  is a positive real parameters.

 For any $i=1,...,N,$ $(-\Delta)^{s}_{a_i(.)}$ are  nonlocal integro-differential operators of elliptic type defined as :
  \begingroup\makeatletter\def\f@size{11}\check@mathfonts
             $$
              \begin{aligned}
              (-\Delta)^s_{a_i(.)}u(x)=2\lim\limits_{\varepsilon\searrow 0} \int_{\R^N\setminus B_\varepsilon(x)} a_i\left( \dfrac{|u(x)-u(y)|}{|x-y|^{s} }\right)\dfrac{u(x)-u(y)}{|x-y|^{s}} \dfrac{dy}{|x-y|^{N+s}},
              \end{aligned}
               $$\endgroup  
 for all $x\in \R^N$, where for any $i=1,..,N,$  $a_i : \R\longrightarrow \R$ which will be specified later.

The most canonical and important example of nonlocal operator is given by the
fractional Laplacian:
$$(-\Delta)^s u(x)=p.v.\int_{\R^N\setminus B_\varepsilon(x)} \dfrac{(u(x)-u(y))}{|x-y|^{N+2s}}dy,$$
which can be obtained by minimizing the Gagliardo seminorm of the fractional
Sobolev space $H^s$, $s\in(0,1)$ (see for instance \cite{11}); from this point of view, it
is important to point out that when $a_i(t)=1$ in \hyperref[P]{$(P_a)$}, then $(-\Delta)^s_{a_i(.)}$ becomes (a multiple of) $(-\Delta)^s$. More generally, when $a_i(t)=t^{p-2}$
in \hyperref[P]{$(P_a)$}, $1<p<\infty$,
we get the eigenvalue problem for the so called fractional $p$-Laplacian
$$(-\Delta)^s_p u(x)=2\lim_{\varepsilon \searrow 0}\int_{\R^N\setminus B_\varepsilon(x)} \dfrac{|u(x)-u(y)|^{p-2}(u(x)-u(y))}{|x-y|^{N+sp}}dy.$$
For the problems involving fractional $p$-laplace operator, we refer the reader to the works \cite{SRB, SR0, SR2}. they use different methods to establish the existence of
  solutions.\\
  
  During the last decade, much attention has been devoted to the study of eigenvalue problems involving non-local operators. The archetypal example being the fractional $p$-Laplacian, has been widely studied and is by now fairly well understood. That is, in the works of Franzina \& Palatucci \cite{fran} and of Lindgren \& Linqvist \cite{lin}. For non-homogeneous problems, eigenvalue problems in classical Orlicz-Sobolev spaces have also been studied thoroughly, see
  for instance, Cl\'ement et al \cite{cl2}, Mih\u{a}ilescu \& R\u{a}dulescu, \cite{ra2}, and Mih\u{a}ilescu et al \cite{ra3}. The theory of eigenvalue problems in fractional Orlicz-Sobolev  spaces is still under development. In this context, Salort \cite{sal2} studied the eigenvalues and minimizers of a fractional nonstandard growth problem in fractional Orlicz-Sobolev spaces and he established several properties on these quantities and their corresponding eigenfunctions. Azroul et al in \cite{srati3} established the existence of two positive eigenvalues for a nonlocal problem, the same authors in \cite{sr5} by means of Ekeland's variational principle and direct variational approach. They obtained the existence of nontrivial
  weak solutions for an eigenvalue problem in fractional Orlicz-Sobolev spaces.

Inspired by these previous works, and by using some different variational methods, the aim of this paper is to study the existence or non-existence of weak solutions for Problem  \hyperref[P]{$(P_a)$} in fractional anisotropic Orlicz-Sobolev spaces. 
 This is the first contribution to studying of nonlocal problems in this class of fractional anisotropic Orlicz Sobolev spaces.\\
           
 This paper is organized as follows : In the second section we  introduce the  fractional anisotropic Orlicz-Sobolev spaces and  we establish some qualitative properties of these
 new spaces. In the third section, we present our main results, and by using some variational methods, we establish the existence or non-existence of eigenvalues for Problem \hyperref[P]{$(P_a)$}.
\section{Variatoinal setting and preliminaries results}

To deal with this situation we define the fractional Orlicz-Sobolev space to investigate Problem \hyperref[P]{$(P_a)$}. Let us recall the definitions and some elementary properties of this spaces. We refer the reader to \cite{1,3,sal1} for further reference and for some of the proofs of the results in this section.\\

 Let $\Omega$ be an open subset of $\R^N$, $N\geqslant 1$. We assume that any $i=1,...,N,$ $a_i : \R\longrightarrow \R$ in \hyperref[P]{$(P_a)$} are such that : $\varphi_i : \R\longrightarrow \R$ defined by : 
  $$
     \varphi_i(t)= \left\{ 
          \begin{array}{clclc}
        a_i(|t|)t   & \text{ for }& t\neq 0, \\\\
          0  & \text{ for } & t=0,
          \end{array}
          \right. 
       $$
is increasing homeomorphisms from $\R$ onto itself. Let 
$$\varPhi_i(t)=\int_{0}^{t}\varphi_i(\tau)d\tau.~~\text{for all}~~ t\in \R,~~i=1,...,N.$$  
Then, for any $i=1,...,N$, $\varPhi_i$, are $N$-functions, see \cite{1}, that is $\varPhi_i :\R^+\longrightarrow \R^+$ are continuous, convex, increasing functions, with $\dfrac{\varPhi_i(t)}{t}\rightarrow 0$ as $t\rightarrow 0$ and $\dfrac{\varPhi_i(t)}{t}\rightarrow \infty$ as $t\rightarrow \infty.$

For the functions $\varPhi_i$ ($i=1,...,N$), introduced above we define the Orlicz spaces :
 $$L_{\varPhi_i} (\Omega)=\left\lbrace u : \Omega \longrightarrow \R \text{ mesurable : } \int_\Omega\varPhi_i(\lambda |u(x)|)dx < \infty \text{ for some } \lambda>0 \right\rbrace. $$
For any $i=1,...,N$, the spaces $L_{\varPhi_i}(\Omega)$ , are  Banach spaces endowed with the Luxemburg norm 
$$||u||_{\varPhi_i}=\inf\left\lbrace \lambda>0 \text{ : }\int_\Omega\varPhi_i\left( \dfrac{|u(x)|}{\lambda}\right) dx\leqslant 1\right\rbrace. $$
 The conjugate $N$-functions of $\varPhi_i$ ($i=1,...,N$), are defined by $\overline{\varPhi_i}(t)=\int_{0}^{t}\overline{\varphi_i}(\tau)d\tau$, where $\overline{\varphi_i} : \R\longrightarrow \R,$ $i=1,...,N$, are given by $\overline{\varphi_i}(t)=\sup\left\lbrace s \text{ : } \varphi_i(s)\leqslant t\right\rbrace.$ Furthermore, it is possible to prove a H\"older type inequality, that is
  \begin{equation}
   \left| \int_{\Omega}uvdx\right| \leqslant 2||u||_{\varPhi_i}||v||_{\overline{\varPhi_i}}\hspace*{0.5cm} \text{  } \forall u \in L_{\varPhi_i}(\Omega)  \text{ and } v\in L_{\overline{\varPhi_i}}(\Omega)~~i=1,...,N.
   \end{equation}
    Throughout this paper, we assume that
 \begin{equation}\label{f1.}
    1<\varphi_i^-:=\inf_{t\geqslant 0}\dfrac{t\varphi_i(t)}{\varPhi_i(t)}\leqslant \varphi_i^+:=\sup_{t\geqslant 0}\dfrac{t\varphi_i(t)}{\varPhi_i(t)}<+\infty~~i=1,...,N. \end{equation}
             The above relation implies that for any $i=1,...,N,$ $\varPhi_i\in \Delta_2$ that is $\varPhi_i$ satisfies the global $\Delta_2$-condition (see \cite{ra}) :
$$
    \varPhi_i(2t)\leqslant K_i\varPhi_i(t)~~ \text{for all}~~ t\geqslant 0~~i=1,...,N,
$$
 where $K_i$ $(i=1,...,N)$, are positive constants. 
 
 Furthermore, for any $i=1,...,N$, we assume that $\varPhi_i$ satisfies the following condition
 \begin{equation}\label{f2.}
 \text{ the functions } [0, \infty) \ni t\mapsto \varPhi_i(\sqrt{t}) \text{ is convex. }
 \end{equation} 
 The  above relation assures that $L_{\varPhi_i}(\Omega)$ is an uniformly convex spaces (see \cite{ra}).
    
\begin{lem}$\label{lemp}$(cf. \cite{ra}).
                        Assume condition $(\ref{f1.})$ is satisfied. Then, for every $b> 1$ and $t\geqslant 0$, we have
                    \begin{equation}\label{3.}
              \varPhi_i(bt)\leqslant b^{\varphi_i^+}\varPhi_i(t),~~ i=1,...,N.
                        \end{equation}
                        \end{lem}
                        
    \begin{defini}               
           Let $A$, $B$ be two N-function. 
                  $A$ is stronger (resp essentially stronger) than $B$,  $A\succ B$ (resp $A\succ\succ B$) in symbols, if 
                $$B(x)\leqslant A(a x ), \text{    } x\geqslant x_0\geqslant 0, $$
                for some (resp for each) $a>0$ and $x_0$ (depending on $a$).
                      \end{defini}
                  \begin{rem}(see. \cite[Section 8.5]{1}).
            $A\succ\succ B$  is equivalent to the condition \\
                           $$\lim_{x\rightarrow \infty}\dfrac{B(\lambda x)}{A(x)}=0,$$
                           for all $\lambda>0$. 
                              \end{rem}

Now, for any $i=1,...,N$, we defined the fractional Orlicz-Sobolev spaces $W^{s}L_{\varPhi_i}(\Omega)$ as follows 
\begingroup\makeatletter\def\f@size{9}\check@mathfonts$$ W^{s}{L_{\varPhi_i}}(\Omega)=\Bigg\{u\in L_{\varPhi_i}(\Omega) ~ :~ \int_{\Omega} \int_{\Omega} \varPhi_i\left( \dfrac{\lambda| u(x)- u(y)|}{|x-y|^{s}}\right) \dfrac{dxdy}{|x-y|^N}< \infty~~ \text{
 for some }\lambda>0 \Bigg\}.
$$\endgroup
This spaces are equipped with the norm,
\begin{equation}\label{6}
||u||_{s,\varPhi_i}=||u||_{\varPhi_i}+[u]_{s,\varPhi_i},
\end{equation}
where $[.]_{s,\varPhi_i}$ is the Gagliardo seminorm, defined by 
$$[u]_{s,\varPhi_i}=\inf \Bigg\{\lambda >0 :  \int_{\Omega} \int_{\Omega} \varPhi_i\left( \dfrac{|u(x)- u(y)|}{\lambda|x-y|^{s}}\right) \dfrac{dxdy}{|x-y|^N}\leqslant 1 \Bigg\}.
$$ 
 By \cite{sal1}, for any $i=1,...,N,$ $W^{s}L_{\varPhi_i}(\Omega)$   are Banach spaces, also  separable (resp. reflexive) spaces if and only if $\varPhi_i\in \Delta_2$ (resp. $\varPhi_i\in \Delta_2$ and $\overline{\varPhi_i}\in \Delta_2$). Furthermore
 if   $\varPhi_i\in \Delta_2$ and $\varPhi_i(\sqrt{t})$ are convex, then  the spaces $W^{s}L_{\varPhi_i}(\Omega)$ are uniformly convex. 
 
 We also define the fractional Orlicz-Sobolev spaces $W^{s}_0L_{\varPhi_i}(\Omega)$, $i=1,...,N$, as the closure of $C^\infty_0(\Omega)$ in $W^{s}L_{\varPhi_i}(\Omega)$. By \cite[Theorem 6]{3}, we can be equivalently renormed by setting $[.]_{s,\varPhi_i}$.

 \begin{prop}\label{pro3} (cf. \cite{3})
             Assume condition $(\ref{f1.})$ is satisfied. Then the following relations holds true,
               {\small\begin{equation*}
              [u]^{\varphi_i^-}_{s,\varPhi_i}\leqslant \phi_i(u) \leqslant [u]^{\varphi_i^+}_{s,\varPhi_i},~~ \forall u\in W^{s}L_{\varPhi_i}(\Omega) \text{ with }[u]_{s,\varPhi_i}>1,~~i=1,...,N,
              \end{equation*} }
             {\small  \begin{equation*}
                 [u]^{\varphi_i^+}_{s,\varPhi_i}\leqslant \phi_i(u) \leqslant [u]^{\varphi_i^-}_{s,\varPhi_i},~~ \forall u\in W^{s}L_{\varPhi_i}(\Omega) \text{ with }[u]_{s,\varPhi_i}<1,~~i=1,...,N,
                 \end{equation*}} 
  where
  $$\phi_i(u)=\int_{\Omega} \int_{\Omega} \varPhi_i\left( \dfrac{\lambda| u(x)- u(y)|}{|x-y|^{s}}\right) \dfrac{dxdy}{|x-y|^N}.$$                 
                              \end{prop}
  \begin{thm} (cf. \cite{3})    \label{3.1..}
             Let $\Omega$ be a bounded open subset of $\R^N$. 
             Then,
             $$C^\infty_0(\Omega)\subset C^2_0(\Omega)\subset W^{s}L_{\varPhi_i}(\Omega)~~i=1,...,N.$$
        \end{thm}
                                      
 In this paper, we assume that there exists $j\in\left\lbrace 1,...,N\right\rbrace $ such that:
      \begin{equation}\label{15}
      \int_{0}^{1} \dfrac{\varPhi_j^{-1}(\tau)}{\tau^{\frac{N+s}{N}}}d\tau<\infty ~~ \text{ and } ~~ \int_{1}^{\infty} \dfrac{\varPhi_j^{-1}(\tau)}{\tau^{\frac{N+s}{N}}}d\tau=\infty.
      \end{equation}    
      We define the inverse Sobolev conjugate $N$-function of $\varPhi_j$ as follows, 
        \begin{equation}\label{c2}
        (\varPhi_j)_*^{-1}(t)=\int_{0}^{t}\dfrac{\varPhi_j^{-1}(\tau)}{\tau^{\frac{N+s}{N}}}d\tau.
        \end{equation}  
        
        \begin{thm}\label{th1}(cf. \cite{3})
         Let $\Omega$  be a bounded open
          subset of  $\R^N$ with $C^{0,1}$-regularity 
            and bounded boundary. If $(\ref{f1.})$ and  $(\ref{15}) $   hold, then 
         \begin{equation}\label{18}
          W^{s}{L_{\varPhi_j}}(\Omega)\hookrightarrow L_{(\varPhi_j)_*}(\Omega).
         \end{equation}
        \end{thm}
        \begin{thm}\label{th2}(cf. \cite{3})
                Let $\Omega$  be a bounded open
                  subset of  $\R^N$ and  $C^{0,1}$-regularity 
                    with bounded boundary. If $(\ref{f1.})$ and  $(\ref{15}) $   hold, then 
                 \begin{equation}\label{27}
                  W^{s}{L_{\varPhi_j}}(\Omega)\hookrightarrow L_{B}(\Omega),
                 \end{equation}
                 is compact for all $B\prec\prec (\varPhi_j)_*$.
                 \end{thm}
                 \begin{rem}
                  We point out certain examples of functions $\varphi : \R \rightarrow \R$ which are odd, increasing homeomorphisms from $\R$ into $\R$ and satisfy conditions $(\ref{f1.})$ and $(\ref{f2.})$ (see \cite{cl2}).
                  
                  $1)$ Let $$\varphi(t)=p|t|^{p-2}t ~~\forall t\in \R,$$
                  with $p>1$, For this function it can be proved that
                  $$\varphi^-=\varphi^+=p.$$
                  Furthermore, in this particular case the corresponding Orlicz space $L_\varPhi(\Omega)$ is the classical  Lebesgue space $L^p(\Omega)$ while the fractional Orlicz-Sobolev space $W^sL_\varPhi(\Omega)$ is the fractional Sobolev spaces $W^{s,p}(\Omega)$.
                  
                  $2)$  Consider 
                  $$\varphi(t)=\log(1+|t|)|t|^{p-2}t, ~~\forall t\in \R$$
                   with $p>1$. In this case it can be proved that 
                   $$\varphi^-=p,~~\varphi^+=p+1.$$
                   
                   $3)$ Let 
                   $$\varphi(t)=\dfrac{|t|^{p-2}t}{\log(1+|t|)}, ~~\text{ if } t \neq 0, \varphi(0)=0$$
                   with $p>2$. In this case we have
                     $$\varphi^-=p-1,~~\varphi^+=p.$$
                  \end{rem} 
                 
Next, we recall some useful properties of variable exponent spaces. For more details we refer the reader to \cite{23,27}, and the references therein.\\ 
  Consider the set
   $$C_+(\overline{\Omega})=\left\lbrace q\in C(\overline{\Omega}): q(x)>1\text{ for all } x \in\overline{\Omega}\right\rbrace .$$
   For all $q\in C_+(\overline{\Omega}) $, we define $$q^{+}= \underset{x\in \overline{\Omega}}{\sup}~q(x) \quad\text{and}\quad q^{-}= \underset{x\in \overline{\Omega}}{\inf}~q(x).$$
For any  $q\in C_+(\overline{\Omega}) $, we define the variable exponent Lebesgue space as $$L^{q(x)}(\Omega)=\bigg\{u:\Omega\longrightarrow \mathbb{R} ~~\text{measurable}: \int_{\Omega}|u(x)|^{q(x)}dx<+\infty
\bigg\}.$$
This vector space endowed with the \textit{Luxemburg norm}, which is defined by
$$\|u\|_{q(x)}= \inf \bigg\{\lambda>0:\int_{\Omega}\bigg|\frac{u(x)}{\lambda}\bigg|^{q(x)}dx \leqslant1 \bigg\}$$
is a separable reflexive Banach space.

 A very important role in manipulating the generalized Lebesgue spaces with variable exponent is played by the modular of the $L^{q(x)}(\Omega)$ space, which defined by
 $$\begin{array}{clc}
 \hspace{-0.8cm}\rho_{q(.)}: L^{q(x)}(\Omega)\longrightarrow\mathbb{R}\\
  \hspace{5.3cm}u\longmapsto\rho_{q(.)}(u)=\displaystyle\int_{\Omega}|u(x)|^{q(x)}dx.
 \end{array}$$
\begin{prop}$\label{anproop5}$
Let $u\in  L^{q(x)}(\Omega) $, then we have
\begin{enumerate}[label=(\roman*)]
\item $\|u\|_{L^{q(x)}(\Omega)}<1$ $(resp. =1, >1)$ $\Leftrightarrow$ $ \rho_{q(.)}(u)<1$ $(resp. =1, >1)$,
\item  $\|u\|_{L^{q(x)}(\Omega)}<1$ $\Rightarrow$ $\|u\|^{q{+}}_{L^{q(x)}(\Omega)}\leqslant \rho_{q(.)}(u)\leqslant \|u\|^{q{-}}_{L^{q(x)}(\Omega)}$,
\item  $\|u\|_{L^{q(x)}(\Omega)}>1$ $\Rightarrow$ $\|u\|^{q{-}}_{L^{q(x)}(\Omega)}\leqslant \rho_{q(.)}(u)\leqslant \|u\|^{q{+}}_{L^{q(x)}(\Omega)}$.
\end{enumerate}
\end{prop}
Now, in order to study Problem \hyperref[P]{$(P_a)$}, it is important to encode the boundary condition $u=0$ in $\R^N\setminus \Omega$  in the weak formulation. In the case of fractional Sobolev space with variable exponent, Azroul et al \cite{SRH} introduced a new function space to study the variational functionals related to the fractional $p(x,.)$-Laplacian operator by observing the interaction between $\Omega$ and $\R^N\setminus \Omega$. Motivated by the above paper, and due to the nonlocality of the operator $(-\Delta)^s_{a_i(.)}$, we introduce the new fractional Orlicz-Sobolev space 
as follows
\begingroup\makeatletter\def\f@size{9}\check@mathfonts $$W^sL_{\varPhi_i}(Q)=\Bigg\{u\in L_{\varPhi_i}(\Omega) ~ :~ \int_{Q}  \varPhi_i\left( \dfrac{\lambda|u(x)- u(y)|}{|x-y|^{s}}\right) \dfrac{dxdy}{|x-y|^N}< \infty~~ \text{
 for some }\lambda>0 \Bigg\},
$$\endgroup
for any $i=1,...,N$, where $Q=\R^{2N}\setminus (C\Omega\times C\Omega)$ with $C\Omega=\R^N \setminus \Omega$. This spaces are equipped with the norm,
\begin{equation}\label{an6}
||u||_{i}=||u||_{\varPhi_i}+[u]_{i},
\end{equation}
where $[.]_{i}$ is the Gagliardo seminorm, defined by 
$$[u]_{i}=\inf \Bigg\{\lambda > 0 :  \int_{Q} \varPhi_i\left( \dfrac{|u(x)- u(y)|}{\lambda|x-y|^{s}}\right) \dfrac{dxdy}{|x-y|^N}\leqslant 1 \Bigg\}.
$$ 
 Similar to the spaces $(W^sL_{\varPhi_i}(\Omega), \|.\|_{s,\varPhi_i})$ we have that $(W^sL_{\varPhi_i}(Q),  \|.\|_{i})$ are a separable reflexive Banach spaces.\\
 
 Now, let $W_0^sL_{\varPhi_i}(Q)$ denotes the following linear subspace of $W^sL_{\varPhi_i}(Q),$
 $$W_0^sL_{\varPhi_i}(Q)=\left\lbrace u\in W^sL_{\varPhi_i}(Q) ~:~ u=0 \text{ a.e in } \R^N \setminus \Omega\right\rbrace $$
 with the norm
 $$[u]_{i}=\inf \Bigg\{\lambda > 0 :  \int_{Q} \varPhi_i\left( \dfrac{|u(x)- u(y)|}{\lambda|x-y|^{s}}\right) \dfrac{dxdy}{|x-y|^N}\leqslant 1 \Bigg\}.
 $$ 
 It is easy to check that $[u]_{i}$ is a norm on $W_0^sL_{\varPhi_i}(Q)$ (see Corollary $\ref{ann}$).
 
 In the following theorem, we compare the spaces $W^sL_{\varPhi_i}(\Omega)$ and $W^sL_{\varPhi_i}(Q)$.
 \begin{thm}\label{an2} For any $i=1,...,N$, the following assertions hold:
 \begin{itemize}
 \item[1)] The continuous embedding $$W^sL_{\varPhi_i}(Q)\subset W^{s}L_{\varPhi_i}(\Omega)$$
  holds true.\\
 \item[2)] If $u\in W_0^sL_{\varPhi_i}(Q)$, then $u\in W^{s}L_{\varPhi_i}(\R^N)$ and 
  $$||u||_{s,\varPhi_i}\leqslant ||u||_{W^sL_{\varPhi_i}(\R^N)}=||u||_{i}.$$
   \end{itemize}
   \end{thm}
   \begin{proof}[\textbf{Proof}]
   $1)$ Let $u\in W^sL_{\varPhi_i}(Q)$, since $\Omega\times \Omega\subsetneq Q,$ then for all $\lambda>0$ we have 
    \begin{equation}\label{an1}
   \int_{\Omega}\int_{\Omega} \varPhi_i\left( \dfrac{|u(x)- u(y)|}{\lambda|x-y|^{s}}\right) \dfrac{dxdy}{|x-y|^N}\leqslant \int_{Q} \varPhi_i\left( \dfrac{|u(x)- u(y)|}{\lambda|x-y|^{s}}\right) \dfrac{dxdy}{|x-y|^N}.    \end{equation}
    We set 
    $$\mathcal{A}^s_{\lambda,\Omega\times \Omega}=\Bigg\{\lambda > 0 :  \int_{\Omega}\int_{\Omega} \varPhi_i\left( \dfrac{|u(x)- u(y)|}{\lambda|x-y|^{s}}\right) \dfrac{dxdy}{|x-y|^N}\leqslant 1 \Bigg\}$$
    and 
     $$\mathcal{A}^s_{\lambda,Q}=\Bigg\{\lambda > 0 :  \int_{Q} \varPhi_i\left( \dfrac{|u(x)- u(y)|}{\lambda|x-y|^{s}}\right) \dfrac{dxdy}{|x-y|^N}\leqslant 1 \Bigg\}.$$
     By $(\ref{an1})$, it is easy to see that $\mathcal{A}^s_{\lambda,Q}\subset \mathcal{A}^s_{\lambda,\Omega\times\Omega}$. Hence 
     \begin{equation}\label{an}
    [u]_{s,\varPhi_i}=\inf\limits_{\lambda>0}\mathcal{A}^s_{\lambda,\Omega\times\Omega}\leqslant[u]_{i}=\inf\limits_{\lambda>0}\mathcal{A}^s_{\lambda,Q}. \end{equation}
  Consequently, by definitions of the norms $\|u\|_{s,\varPhi_i}$ and $\|u\|_{i},$ we obtain
    $$ \|u\|_{s,\varPhi_i}\leqslant \|u\|_{i}<\infty.$$
   $2)$ Let $u\in W_0^sL_{\varPhi_i}(Q)$, then $u=0$ in $\R^N\setminus \Omega$. So, $\|u\|_{L_{\varPhi_i}(\Omega)}=\|u\|_{L_{\varPhi_i}(\R^N)}.$ Since 
    $$\int_{\R^{2N}} \varPhi_i\left( \dfrac{|u(x)- u(y)|}{\lambda|x-y|^{s}}\right) \dfrac{dxdy}{|x-y|^N}=\int_{Q} \varPhi_i\left( \dfrac{|u(x)- u(y)|}{\lambda|x-y|^{s}}\right) \dfrac{dxdy}{|x-y|^N}$$
    for all $\lambda>0$. Then $[u]_{W^sL_{\varPhi_i}(\R^N)}=[u]_{i}$. Thus, we get
    
     $$||u||_{s,\varPhi_i}\leqslant ||u||_{W^sL_{\varPhi_i}(\R^N)}=||u||_{i}.$$
    \end{proof}
    \begin{coro}\label{ann}(Poincar\'{e} inequality)
    Let $\Omega$ be a bounded subset of $\R^N$. Then there exists a positive
    constant $c$ such that,
    $$
    \|u\|_{\varPhi_i}\leqslant c[u]_i, ~~~~\forall u\in W^s_0L_{\varPhi_i}(Q).$$
    \end{coro}
     \begin{proof}[\textbf{Proof}]
 Let $u\in W^s_0L_{\varPhi_i}(Q)$, by Theorem $\ref{an2}$, we have $u\in W^s_0L_{\varPhi_i}(\Omega)$. Then by \cite[Theorem 6]{3},  there exists a positive
     constant $c$ such that,
     $$
     \|u\|_{\varPhi_i}\leqslant c[u]_{s,\varPhi_i}.$$
    Combining the above inequality with $(\ref{an})$, we obtain that 
     $$
         \|u\|_{\varPhi_i}\leqslant c[u]_i, ~~~~\forall u\in W^s_0L_{\varPhi_i}(Q).$$
     \end{proof}    
 Now, we introduce a natural  fractional anisotropic Orlicz-Sobolev spaces $W_0^sL_{\overrightarrow{\varPhi}}(Q)$, that will enable us to study Problem \hyperref[P]{$(P_a)$}. For this purpose, let us denote by $\overrightarrow{\varPhi} : \Omega \longrightarrow \R^N$ the vectorial function $\overrightarrow{\varPhi}=(\varPhi_1,...,\varPhi_N)$.
 We define $W_0^sL_{\overrightarrow{\varPhi}}(Q)$, the fractional anisotropic Orlicz-Sobolev space as follows the closure of $C^\infty_0(\Omega)$ with respect to the norm:
 $$\|u\|_{\overrightarrow{\varPhi}}=\sum_{i=1}^{N}[u]_{i}.$$
Denoting $X=L_{\varPhi_1}(\Omega)\times...\times L_{\varPhi_N}(\Omega)$ and considering  the operator $T : W_0^sL_{\overrightarrow{\varPhi}}(Q)\longrightarrow X,$ defined by $T(u)=\left( D^su,...,D^su\right)$
where $$D^su=\dfrac{u(x)-u(y)}{|x-y|}.$$ It is clear that $W_0^sL_{\overrightarrow{\varPhi}}(Q)$ and $X$ are isometric by $T$, since $$T(u)=\sum_{i=1}^{N}[u]_{i}=\|u\|_{\overrightarrow{\varPhi}}.$$
Thus $T(W_0^sL_{\overrightarrow{\varPhi}}(Q))$ is a closed subspace of $X$, which is a reflexive Banach space. By  \cite[Proposition III.17]{breziz}, it follows that  $T(W_0^sL_{\overrightarrow{\varPhi}}(Q))$ is reflexive and consequently  $W_0^sL_{\overrightarrow{\varPhi}}(Q)$ is also reflexive Banach space.
 
On the other hand, in order to facilitate the manipulation of the space $W_0^sL_{\overrightarrow{\varPhi}}(Q)$, we introduce $\overrightarrow{\varphi}^+, \overrightarrow{\varphi}^-\in \R^N$ as
$$ \overrightarrow{\varphi}^+=(\varphi_1^+,...,\varphi_N^+),~~ \text{and}~~ \overrightarrow{\varphi}^-=(\varphi_1^-,...,\varphi_N^-),$$
 and $\varphi^+_{max}$, $\varphi^-_{max}$, $\varphi^-_{min}$ as 
 $$ \varphi^+_{max}=\max\left\lbrace \varphi_1^+,...,\varphi_N^+\right\rbrace,~~ \varphi^-_{max}=\max\left\lbrace \varphi_1^-,...,\varphi_N^-\right\rbrace,$$ $$ \varphi^-_{min}=\min\left\lbrace \varphi_1^-,...,\varphi_N^-\right\rbrace.$$
 Throughout this paper we assume that 
 \begin{equation}\label{an8}
   \lim\limits_{t\rightarrow \infty}\dfrac{|t|^{q^+}}{(\varPhi_j)_*(kt)}=0 ~~\forall k>0,
   \end{equation}
   where $j\in\left\lbrace 1,...,N\right\rbrace $ is given in $(\ref{15})$.
 \begin{thm}\label{anth5}
   Let $\Omega$  be a bounded open
             subset of  $\R^N$ with $C^{0,1}$-regularity 
               and bounded boundary. Then 
             the embedding    
                    \begin{equation}\label{an27}
                     W_0^{s}{L_{\overrightarrow{\varPhi}}}(Q)\hookrightarrow L^{q(x)}(\Omega)
                    \end{equation}
                    is compact.
                    \end{thm}
\begin{proof}[\textbf{Proof}]
Let $u\in  W_0^{s}{L_{\overrightarrow{\varPhi}}}(Q)$, so $u\in  W_0^{s}{L_{\varPhi_i}}(Q)$ for any $i=1,...N$, and by Theorem $\ref{an2}$ we have $u\in  W_0^{s}{L_{\varPhi_i}}(\Omega)$ for any $i=1,...N$, then for $j\in \left\lbrace 1,...,N\right\rbrace $ given by ($\ref{an8}$), we can apply Theorem $\ref{th2}$, and we have
$$\|u\|_{q^+}\leqslant c[u]_{s,\varPhi_j}\leqslant c[u]_{j}\leqslant c\sum_{i=1}^{N}[u]_{i}=c\|u\|_{\overrightarrow{\varPhi}}.$$
This implies that 
$$ W_0^{s}{L_{\overrightarrow{\varPhi}}}(Q)\hookrightarrow L^{q^+}(\Omega).$$
That fact combined with the continuous embedding of $L^{q^+}(\Omega)$ in $L^{q(x)}(\Omega)$  ensures that $W_0^{s}L_{\overrightarrow{\varPhi}}(Q)$ is compactly embedded in $L^{q(x)}(\Omega)$.
\end{proof}
 We put, $$\Psi(u)=\displaystyle\int_{Q}\sum_{i=1}^{N}\varPhi_i\left(\dfrac{|u(x)- u(y)|}{|x-y|^{s}}\right) \dfrac{dxdy}{|x-y|^N}.$$ 
 Then similar to proof of Proposition $\ref{pro3}$, we have the following proposition.   
 \begin{prop}\label{anpro3}
             Assume condition $(\ref{f1.})$ is satisfied. Then the following relations holds true,
              \begin{equation}
              \|u\|^{\varphi^-_{\min}}_{\overrightarrow{\varPhi}}\leqslant \Psi(u) \leqslant \|u\|^{\varphi^+_{\max}}_{\overrightarrow{\varPhi}} \text{   ,  } \forall u\in W_0^{s}L_{\overrightarrow{\varPhi}}(Q) \text{ with }\|u\|_{\overrightarrow{\varPhi}}>1,
              \end{equation} 
               \begin{equation}
                 \|u\|^{\varphi^+_{\max}}_{\overrightarrow{\varPhi}}\leqslant \Psi(u) \leqslant \|u\|^{\varphi^-_{\min}}_{\overrightarrow{\varPhi}} \text{   ,  } \forall u\in W_0^{s}L_{\overrightarrow{\varPhi}}(Q) \text{ with }\|u\|_{\overrightarrow{\varPhi}}<1.
                 \end{equation}            
                              \end{prop}
   Finally, the proof of our main results is based on the following  mountain pass theorem, Ekeland's variational principle theorem and direct variational approach.
   \begin{thm}\cite{110}\label{an2.2}
         Let $X$ be a real Banach space and $I \in C^1(X,\R)$ with $I(0)=0$. Suppose that the following conditions hold:
         
         $(G_1)$ \label{G1} There exist
         $\rho>0 \text{  and } r>0 \text{ such that } I(u)\geqslant r \text{ for } ||u||=\rho$.
          
         $(G_2)$ \label{G2} There exists
         $e \in X \text{ with } ||e||>\rho \text{  such that } I(e)\leqslant 0$.\\
         Let
         $$c:=\inf_{\gamma \in \Gamma} \max_{t\in [0,1]}I (\gamma (t)) \text{ with } \Gamma=\left\lbrace \gamma \in C([0,1],X); \gamma(0)=0 , \gamma(1)=e \right\rbrace.$$
         Then there exists a sequence $\left\lbrace u_n\right\rbrace $ in $X$ such that 
         $$I(u_n)\rightarrow c \text{ \hspace{0.2cm} \text{and} \hspace{0.2cm}} I'(u_n)\rightarrow 0.$$
        
         \end{thm}
    \begin{thm}\label{anek}(see : \cite{ek})
    Let V be a complete metric space and $F : V \longrightarrow \R\cup \left\lbrace +\infty\right\rbrace$ be a lower semicontinuous functional on $V$, that is bounded below and not identically equal to $+\infty$. Fix $\varepsilon>0$ and a  point $u\in V$ 
      such that
     $$F(u)\leqslant \varepsilon +\inf\limits_{x\in V}F(x).$$ Then for every $\gamma > 0$,
      there exists some point $v\in V$ such that :
      $$F(v)\leqslant F(u),$$
      $$d(u,v)\leqslant \gamma$$
      and for all $w\neq v$
      $$F(w)> F(v)-\dfrac{\varepsilon}{\gamma}d(v,w).$$
    \end{thm}
    
     \begin{thm}\label{andr} (see :\cite{110})
         Suppose that $X$ is a reflexive Banach space with norm $||.||$ and let
         $V\subset X$ be a weakly closed subset of $X$. Suppose $E : V \longrightarrow \R \cup \left\lbrace +\infty\right\rbrace $ is coercive
         and (sequentially) weakly lower semi-continuous on $V$ with respect to $X$, that
         is, suppose the following conditions are fulfilled:
      \begin{itemize}
         \item[$\bullet$] $E(u)\rightarrow \infty$ as $||u||\rightarrow \infty$, $u\in V$.
         
         \item[$\bullet$]  For any $u\in V$, any sequence $\left\lbrace u_n\right\rbrace $ in $V$ such that $u_n\rightharpoonup u$ weakly in $X$
         there holds:
        $$E(u)\leqslant \liminf_{n\rightarrow \infty}E(u_n).$$
          \end{itemize} 
         Then $E$ is bounded from below on $V$ and attains its infimum  in $V$.     
          \end{thm}
         
   \section{ Main results and proofs}
   To simplify the notation, we ask
         $$D^{s}u:=\dfrac{u(x)-u(y)}{|x-y|^{s}}~~ \text{ and }~~ d\mu=\dfrac{dxdy}{|x-y|^N} ~~ \forall (x,y)\in Q. $$ 
         
 The dual space of $\left(W_0^{s}L_{\overrightarrow{\varPhi}}(Q),||.||_{\overrightarrow{\varPhi}}\right) $  is denoted by $\left((W_0^{s}L_{\overrightarrow{\varPhi}}(Q))^*,||.||_{\overrightarrow{\varPhi},*}\right) $.

 \begin{defini}
 We say that $\lambda\in \R$ is an eigenvalue of Problem \hyperref[P]{$(P_a)$} if there exists $u\in W_0^{s}L_{\overrightarrow{\varPhi}}(Q)\setminus \left\lbrace 0\right\rbrace$ such that 
 $$\int_{Q} \sum_{i=1}^{N}a_i(|D^{s}u|)  D^{s}u D^{s}vd\mu-\lambda\int_{\Omega}|u|^{q(x)-2}uvdx=0$$
 for all $v\in W_0^{s}L_{\overrightarrow{\varPhi}}(Q)$.
 \end{defini}
 
  We point that if $\lambda$ is an eigenvalue of Problem \hyperref[P]{$(P_a)$} then the corresponding $u\in W_0^{s}L_{\overrightarrow{\varPhi}}(Q)\setminus\left\lbrace 0\right\rbrace $ is a weak solution of \hyperref[P]{$(P_a)$}.\\

  The main results in this paper are given by the following theorems.
  \begin{thm}\label{anth1}
  Assume that the function $q\in C(\overline{\Omega})$ verifies the hypothesis
  \begin{equation}
  \varphi^+_{\max}<q^-.
  \end{equation}
  Then for all $\lambda>0$ is an eigenvalue of Problem \hyperref[P]{$(P_a)$}.
  \end{thm}  
   \begin{thm}\label{anth2}
    Assume that the function $q\in C(\overline{\Omega})$ verifies the hypothesis
    \begin{equation}\label{an15}
    q^-<\varphi^-_{\min}.
    \end{equation}
    Then there exists $\lambda_*>0$ such that for any $\lambda\in (0,\lambda_*)$ is an eigenvalue of Problem \hyperref[P]{$(P_a)$}.
    \end{thm}
    \begin{thm}\label{anth3}
      Assume that the function $q\in C(\overline{\Omega})$ satisfies the conditions 
      \begin{equation}
     q^+<\varphi^-_{\min}.  \end{equation}
      Then there exist two positive constants $\lambda_*>$ and $\lambda^*>0$ such that $\lambda\in (0,\lambda_*)\cup[\lambda^*,\infty)$ is an eigenvalue of Problem \hyperref[P]{$(P_a)$}.
      \end{thm}
     In order to state the next result, we define
      \begin{equation}
      \lambda_1=\inf\limits_{u\in W_0^{s}L_{\overrightarrow{\varPhi}}(Q)\setminus\left\lbrace 0\right\rbrace }\dfrac{ \displaystyle\int_{Q}\sum_{i=1}^{N}\varPhi_i(|D^{s}u|)d\mu} {\displaystyle\int_{\Omega}\dfrac{1}{q(x)}|u|^{q(x)}dx}
      \end{equation}
      and 
      \begin{equation}
            \lambda_0=\inf\limits_{u\in W_0^{s}L_{\overrightarrow{\varPhi}}(Q)\setminus\left\lbrace 0\right\rbrace }\dfrac{\displaystyle\int_{Q} \sum_{i=1}^{N}a_i(|D^{s}u|)|D^{s}u|^2d\mu} {\displaystyle\int_{\Omega}|u|^{q(x)}dx}.
            \end{equation}
  \begin{thm}\label{anth4}
  Assume that there exist $j_1, j_2, k\in \left\lbrace 1,...,N\right\rbrace$ such that
  \begin{equation}\label{an21}
 \varphi_{j_1}^-=q^-,~~\text{} ~~ \varphi^+_{j_2}=q^+,
  \end{equation}
  and
  \begin{equation}\label{an7}
  q^+<\varphi^-_k.
  \end{equation}
  Then $0<\lambda_0\leqslant\lambda_1$, and  any $\lambda\in [ \lambda_1,\infty)$ is an eigenvalue of Problem \hyperref[P]{$(P_a)$} while no $\lambda\in (0,\lambda_0)$ can be an eigenvalue of Problem \hyperref[P]{$(P_a)$}.
  \end{thm}
  \begin{rem}\label{rem}  
We are not able to tell whether $\lambda_0=\lambda_1$ or $\lambda_0<\lambda_1$. Then, an open problem worth proposing is the study of the existence of eigenvalues of the Problem \hyperref[P]{$(P_a)$} in the interval $[\lambda_0,\lambda_1]$.
 \end{rem}
     \subsection{Auxiliary results}
  In order to prove our main results, we introduce the following functionals $J, I, J_1, I_1 : W^{s}_0L_{\overrightarrow{\varPhi}}(Q) \longrightarrow \R$ by

  $J(u)= \displaystyle\int_{Q}\sum_{i=1}^{N}\varPhi_i\left(|D^su| \right)d\mu$

$I(u)=\displaystyle\int_{\Omega}\dfrac{1}{q(x)}|u|^{q(x)}dx$

$J_1(u)\displaystyle\int_{Q} \sum_{i=1}^{N}a_i(|D^{s}u|)|D^{s}u|^2d\mu$

$I_1(u)=\displaystyle\int_{\Omega}|u|^{q(x)}dx$.

By a standard argument to  \cite{sr5} and \cite{3}, we have $J,I\in C^1(W^{s}_0L_{\overrightarrow{\varPhi}}(Q),\R)$,
$$\left\langle J'(u),v\right\rangle =\int_{Q} \sum_{i=1}^{N}a_i(|D^{s}u|)  D^{s}u D^{s}vd\mu$$
 and 
$$\left\langle I'(u),v\right\rangle =\displaystyle\int_{\Omega}|u|^{q(x)-2}uvdx,$$
for all $u,v\in W^{s}_0L_{\overrightarrow{\varPhi}}(Q)$.\\
Next, for each $\lambda\in \R$, we define the energetic function associated with Problem \hyperref[P]{$(P_a)$}, 
 $T_\lambda : W^{s}_0L_{\overrightarrow{\varPhi}}(Q)\longrightarrow \R $ by 
   $$T_\lambda(u)=J(u)-\lambda I(u).$$
  Clearly, $T_\lambda \in C^1(W^{s}_0L_{\overrightarrow{\varPhi}}(Q), \R)$ and 
  $$\left\langle T'_\lambda (u),v\right\rangle = \left\langle J' (u),v\right\rangle-\lambda\left\langle I' (u),v\right\rangle$$
  for all $u,v\in W^{s}_0L_{\overrightarrow{\varPhi}}(Q)$. 
  \begin{lem}\label{anlem1}
  Assume that the hypothesis of Theorem $\ref{anth1}$ is fulfilled. Then, there exist $\eta>0$ and $\alpha>0$, such that $T_\lambda(u)\geqslant \alpha>0$ for any $u\in W^{s}_0L_{\overrightarrow{\varPhi}}(Q)$ with $\|u\|_{\overrightarrow{\varPhi}}=\eta$.
 \end{lem}
 \begin{proof}[\textbf{Proof}]
 First, we point out that
 \begin{equation*}
 |u(x)|^{q^-}+|u(x)|^{q^+}\geqslant |u(x)|^{q(x)}~~\forall x\in \overline{\Omega}.
 \end{equation*}
 Using the above inequality and the definition of $T_\lambda$, we find that 
 \begin{equation}\label{an3}
 T_\lambda(u)\geqslant \int_{Q}\sum_{i=1}^{N}\varPhi_i\left(|D^su| \right)d\mu-\dfrac{\lambda}{q^-}\left( \|u\|_{q^-}^{q^-}+\|u\|_{q^+}^{q^+}\right) 
 \end{equation}
 for any $u\in W^{s}_0L_{\overrightarrow{\varPhi}}(Q)$. Since $W^{s}_0L_{\overrightarrow{\varPhi}}(Q)$ is continuously embedded in $L^{q^{+}}(\Omega)$ and in $L^{q^{-}}(\Omega)$ it follows that there exist two positive constants $c_1$ and $c_2$ such that
 \begin{equation}\label{an24}
 ||u||_{\overrightarrow{\varPhi}}\geqslant c_1||u||_{q^+} ~~\forall u\in W^{s}_0L_{\overrightarrow{\varPhi}}(Q)\end{equation}
 and 
 \begin{equation}\label{an25}
 ||u||_{\overrightarrow{\varPhi}}\geqslant c_2||u||_{q^-} ~~\forall u\in W^{s}_0L_{\overrightarrow{\varPhi}}(Q).\end{equation}
 Next, for $u\in W^{s}_0L_{\overrightarrow{\varPhi}}(Q)$ with $\|u\|_{\overrightarrow{\varPhi}}<1$, so, we have $[u]_{i}<1$ for any 
 $i=1,...,N$. Then
 \begin{equation}\label{an4}
 \begin{aligned}
 \dfrac{\|u\|^{\varphi^+_{\max}}_{\overrightarrow{\varPhi}}}{N^{\varphi_{\max}^+-1}} & =N\left( \sum_{i=1}^{N}\dfrac{1}{N}[u]_{i}\right)^{\varphi^+_{\max}}\\
 & \leqslant \sum_{i=1}^{N}[u]_{i}^{\varphi^+_{\max}}\\
  & \leqslant \sum_{i=1}^{N}[u]_{i}^{\varphi^+_{i}}\\
  &\leqslant \int_{Q}\sum_{i=1}^{N}\varPhi_i\left(|D^su| \right)d\mu.
  \end{aligned}
 \end{equation}
 Relations $(\ref{an3})-(\ref{an4})$, imply that
$$ 
\begin{aligned}
T_\lambda(u)& \geqslant \dfrac{\|u\|^{\varphi^+_{\max}}_{\overrightarrow{\varPhi}}}{N^{\varphi_{\max}^+-1}}-\dfrac{\lambda}{q^-}\left( (c_1\|u\|_{\overrightarrow{\varPhi}})^{q^-}+(c_2\|u\|_{\overrightarrow{\varPhi}})^{q^+}\right)\\
& = \left( c_3-c_4 \|u\|^{q^+-\varphi^+_{\max}}_{\overrightarrow{\varPhi}}-c_5 \|u\|^{q^--\varphi^+_{\max}}_{\overrightarrow{\varPhi}}\right) \|u\|^{\varphi^+_{\max}}_{\overrightarrow{\varPhi}} 
\end{aligned}
$$
for any $u\in W^{s}_0L_{\overrightarrow{\varPhi}}(Q)$ with $\|u\|_{\overrightarrow{\varPhi}}<1$, where $c_3, c_4$, and $c_5$ are positive constants. Since the function $g : [0,1]\longrightarrow \R$ defined by:
$$g(t)=c_3-c_4 t^{q^+-\varphi^+_{\max}}-c_5 t^{q^--\varphi^+_{\max}}$$
is positive in a neighborhood of the origin, the conclusion of the lemma follows at once. 
\end{proof}
  \begin{lem}\label{anlem2}
    Assume that the hypothesis of Theorem $\ref{anth1}$ is fulfilled. Then, there exist $e>0$ with $\|e\|_{\overrightarrow{\varPhi}}=\eta$ (where $\eta$ is given by Lemma $\ref{anlem1}$) such that $T_\lambda(e)<0$.
   \end{lem}
   \begin{proof}[\textbf{Proof}]
   Let $\theta\in C_0^\infty(\Omega)$, $\theta\geqslant 0$ and $\theta\neq 0$, be fixed, and let $t>1$. Using Lemma $\ref{lemp}$, we find that 
   $$
   \begin{aligned}
  T_\lambda(t\theta) &=\displaystyle\int_{Q}\sum_{i=1}^{N}\varPhi_i\left(t|D^s\theta| \right)d\mu-\lambda\displaystyle\int_{\Omega}\dfrac{1}{q(x)}t^{q(x)}|\theta|^{q(x)}dx\\
  &\leqslant \displaystyle\int_{Q}\sum_{i=1}^{N}t^{\varphi_i^+}\varPhi_i\left(|D^s\theta| \right)d\mu-\lambda\displaystyle\int_{\Omega}\dfrac{t^{q^-}}{q^+}|\theta|^{q(x)}dx\\
    &\leqslant t^{\varphi_{\max}^+}\displaystyle\int_{Q}\sum_{i=1}^{N}\varPhi_i\left(|D^s\theta| \right)d\mu-\dfrac{\lambda t^{q^-}}{q^+}\displaystyle\int_{\Omega}|\theta|^{q(x)}dx.
    \end{aligned}
   $$
   Since $q^->\varphi_{\max}^+$, it is clear that $\lim\limits_{t\rightarrow \infty}T_\lambda(t\theta)=-\infty$. Then, for $t>1$ large enough, we can take $e=t \theta$ such that $\|e\|_{\overrightarrow{\varPhi}}=\eta$ and  $T_\lambda(e)<0$. This completes the proof.
      \end{proof}
  \begin{lem}\label{anlem3}
Assume that the hypothesis of Theorem $\ref{anth2}$ is fulfilled. Then, there exists $\lambda_*>0$ such that for any
$\lambda\in (0,\lambda_*)$, there are $\rho, \alpha>0$, such that $T_\lambda(u)\geqslant \alpha>0$ for any $u\in W^{s}_0L_{\overrightarrow{\varPhi}}(Q)$ with $||u||_{\overrightarrow{\varPhi}}=\rho$.
  \end{lem}  
  \begin{proof}[\textbf{Proof}]
  Since $W^{s}_0L_{\overrightarrow{\varPhi}}(Q)$ is continuously embedded in $L^{q(x)}(\Omega)$, it follows that there exists a positive constant $c_1$  such that
   \begin{equation}\label{an11}
  ||u||_{\overrightarrow{\varPhi}}\geqslant c_1||u||_{q(x)} ~~\forall u\in W^{s}_0L_{\overrightarrow{\varPhi}}(Q)\end{equation} 
   we fix $\rho \in (0,1)$ such that $\rho<\dfrac{1}{c_1}$. Then relation $(\ref{an11})$ implies that 
   $$\|u\|_{q(x)}<1~~\text{for all } u\in W^{s}_0L_{\overrightarrow{\varPhi}}(Q) ~~\text{with } ||u||_{\overrightarrow{\varPhi}}=\rho.$$
   Then, we can apply Proposition $\ref{anproop5}$, and we have
   \begin{equation}\label{an12}
   \int_\Omega |u(x)|^{q(x)}dx\leqslant \|u\|_{q(x)}^{q^-}~~\text{ for all } u\in W^{s}_0L_{\overrightarrow{\varPhi}}(Q)~~ \text{ with } ||u||_{\overrightarrow{\varPhi}}=\rho.
   \end{equation}
   Relation $(\ref{an11})$  and $(\ref{an12})$ implies that
    \begin{equation}\label{an13}
      \int_\Omega |u(x)|^{q(x)}dx\leqslant c_1^{q^-}\|u\|_{\overrightarrow{\varPhi}}^{q^-}~~\text{ for all } u\in W^{s}_0L_{\overrightarrow{\varPhi}}(Q)~~ \text{ with } ||u||_{\overrightarrow{\varPhi}}=\rho.
      \end{equation}
   Taking into account Relations $(\ref{an4})$ and $(\ref{an13})$, we deduce that for any $u \in W^{s}_0L_{\overrightarrow{\varPhi}}(Q)$ with  $||u||_{\overrightarrow{\varPhi}}=\rho$, the following inequalities hold true:
   $$
   \begin{aligned}
   T_\lambda(u) &\geqslant  \dfrac{\|u\|^{\varphi^+_{\max}}_{\overrightarrow{\varPhi}}}{N^{\varphi_{\max}^+-1}}-\dfrac{\lambda}{q^-}\int_\Omega|u(x)|^{q(x)}dx\\
   & \geqslant \dfrac{\|u\|^{\varphi^+_{\max}}_{\overrightarrow{\varPhi}}}{N^{\varphi_{\max}^+-1}}-\dfrac{\lambda c_1^{q-}}{q^-}\|u\|^{q^-}_{\overrightarrow{\varPhi}}\\
  & =\rho^{q^-}\left(\dfrac{\rho^{\varphi^+_{\max}-q^-}}{N^{\varphi_{\max}^+-1}}- \dfrac{\lambda c_1^{q-}}{q^-}\right).
   \end{aligned}
   $$
   Hence, if we define
  \begin{equation}\label{an40}
  \lambda_*=\dfrac{\rho^{\varphi^+_{\max}-q^-}}{2c_1^{q^-} N^{\varphi_{\max}^+-1}}q^-.\end{equation}
  Then, for any $\lambda\in (0,\lambda_*)$ and $u\in W^{s}_0L_{\overrightarrow{\varPhi}}(Q)$ with  $||u||_{\overrightarrow{\varPhi}}=\rho$, we have
  $$T_\lambda(u)\geqslant \alpha>0,$$
  such that 
  $$\alpha=\dfrac{\rho^{\varphi^+_{\max}}}{2 N^{\varphi^+_{\max}-1}}.$$
  This completes the proof.
     \end{proof} 
   \begin{lem}\label{anlem4}
   Assume that the hypothesis of Theorem $\ref{anth2}$ is fulfilled. Then, there exists $\phi>0$ such that $\phi\geqslant 0$,  $\phi\neq 0$, and $T_\lambda(t\phi)<0$ for $t>0$ small enough.
     \end{lem}    
   \begin{proof}[\textbf{Proof}] 
 By assumption   $(\ref{an15})$ we can chose $\varepsilon_0>0$ such that $q^-+\varepsilon_0<\varphi^-_{\min}$. On the other hand, since $q\in C(\overline{\Omega})$, it follows that there exists an open set $\Omega_0\subset \Omega$ such that $|q(x)-q^-|<\varepsilon_0$ for all $x\in \Omega_0$. Thus, $q(x)\leqslant q^-+\varepsilon_0<\varphi^-_{\min}$ for all $x\in \Omega_0$. 
 Let $\phi\in C_0^\infty(\Omega)$ be such that $supp(\phi)\supset \overline{\Omega_0}$, $\phi(x)=1$ for all $x\in \overline{\Omega_0}$, and $0\leqslant \phi\leqslant 1$ in $\overline{\Omega_0}$. Then, using the above information and the definition of $\varphi_i^-$, for any $t\in(0,1)$, we have
 $$
 \begin{aligned}
 T_\lambda(t\phi)& =\displaystyle\int_{Q}\sum_{i=1}^{N}\varPhi_i\left(t|D^s\phi| \right)d\mu-\lambda\displaystyle\int_{\Omega}\dfrac{1}{q(x)}t^{q(x)}|\phi|^{q(x)}dx\\
 &\leqslant \displaystyle\int_{Q}\sum_{i=1}^{N}t^{\varphi_i^-}\varPhi_i\left(|D^s\phi| \right)d\mu-\lambda\displaystyle\int_{\Omega_0}\dfrac{t^{q(x)}}{q(x)}|\phi|^{q(x)}dx\\
     &\leqslant t^{\varphi_{\min}^-}\displaystyle\int_{Q}\sum_{i=1}^{N}\varPhi_i\left(|D^s\phi| \right)d\mu-\dfrac{\lambda t^{q^-+\varepsilon_0}}{q^+}\displaystyle\int_{\Omega_0}|\phi|^{q(x)}dx.
     \end{aligned}
    $$
    Therefore $T_\lambda(t\phi)<0$, for $t<\delta^{1/(\varphi^-_{\min}-q^--\varepsilon_0)}$ with
    $$0<\delta<\min\left\lbrace 1,~~ \dfrac{\dfrac{\lambda}{q^+}\displaystyle\int_{\Omega_0}|\phi|^{q(x)}dx}{\displaystyle\int_{Q}\sum_{i=1}^{N}\varPhi_i\left(|D^s\phi| \right)d\mu}\right\rbrace. $$
   This is possible, since we claim that 
   $$  \displaystyle\int_{Q}\sum_{i=1}^{N}\varPhi_i\left(|D^s\phi| \right)d\mu>0.$$
   Indeed, it is clear that
   $$\int_{\Omega_0}|\phi|^{q(x)}dx\leqslant \int_{\Omega}|\phi|^{q(x)}dx\leqslant \int_{\Omega}|\phi|^{q^-}dx.$$
   On the other hand, since  $W^{s}_0L_{\overrightarrow{\varPhi}}(Q)$ is continuously embedded in $L^{q^-}(\Omega)$, it follows that there exists a positive constant $c$  such that 
   $$|\phi|_{q^-}\leqslant c ||\phi||_{\overrightarrow{\varPhi}}.$$
   The last two inequalities imply that
   $$\|\phi\|_{\overrightarrow{\varPhi}}>0$$
   and combining this fact with Proposition $\ref{anpro3}$, the claim follows at once. The proof of the lemma is now completed.
   \end{proof}  
 \begin{lem}\label{anlem5}
    Assume that the hypothesis of Theorem $\ref{anth3}$ is fulfilled. Then, the function $T_\lambda$ is coercive on $W^{s}_0L_{\overrightarrow{\varPhi}}(Q)$. 
      \end{lem}  
         \begin{proof}[\textbf{Proof}]  
  By relations $(\ref{an24})$ and $(\ref{an25})$, we deduce that, for all $u\in W^{s}_0L_{\overrightarrow{\varPhi}}(Q)$,
  \begin{equation}\label{an20}
  T_\lambda(u)\geqslant \int_{Q}\sum_{i=1}^{N}\varPhi_i\left(|D^su| \right)d\mu-\dfrac{\lambda}{q^-}\left[ (c_1\|u\|_{\overrightarrow{\varPhi}})^{q^-}+(c_2\|u\|_{\overrightarrow{\varPhi}})^{q^+}\right].
  \end{equation}    
Now, we focus our attention on the elements $u\in W^{s}_0L_{\overrightarrow{\varPhi}}(Q)$ with $\|u\|_{\overrightarrow{\varPhi}}> 1$. Denoting 
$$
     \alpha_i= \left\{ 
          \begin{array}{clclc}
        \varphi^+_{\max}   & \text{ if }& [u]_{i}<1, \\\\
          \varphi^-_{\min}  & \text{ if } & [u]_{i}>1.
          \end{array}
          \right.
       $$
Then, we have
  \begin{equation}\label{-}
    \begin{aligned}
    \int_{Q}\sum_{i=1}^{N}\varPhi_i\left(|D^su| \right)d\mu &\geqslant \sum_{i=1}^{N} [u]_{i}^{\alpha_i}\\
    & \geqslant \sum_{i=1}^{N} [u]_{i}^{\varphi^-_{\min}}-\sum_{i=1}^{N}\left(  [u]_{i}^{\varphi^-_{\min}}-[u]_{i}^{\varphi^+_{\max}}\right)\\
    & \geqslant \dfrac{\|u\|^{\varphi^-_{\min}}_{\overrightarrow{\varPhi}}}{N^{\varphi_{\min}^--1}}-N.
    \end{aligned}
      \end{equation}
  Therefore, the above inequality with $(\ref{an20})$, implies that
    \begin{equation}
   T_\lambda(u)\geqslant \dfrac{\|u\|^{\varphi^-_{\min}}_{\overrightarrow{\varPhi}}}{N^{\varphi_{\min}^--1}}-N-\dfrac{\lambda}{q^-}\left[ (c_1\|u\|_{\overrightarrow{\varPhi}})^{q^-}+(c_2\|u\|_{\overrightarrow{\varPhi}})^{q^+}\right] 
    \end{equation}    
  for $u\in W^{s}_0L_{\overrightarrow{\varPhi}}(Q)$ with $\|u\|_{\overrightarrow{\varPhi}}> 1$. Since $\varphi^-_{\min}>q^+\geqslant q^-$, we infer $T_\lambda(u)\rightarrow \infty$ as $\|u\|_{\overrightarrow{\varPhi}}\rightarrow \infty$. Then $T_\lambda$ is coercive on  $ W^{s}_0L_{\overrightarrow{\varPhi}}(Q)$.   
            \end{proof} 
  \begin{lem}\label{anlem6}
     Assume that condition $(\ref{an21})$  is fulfilled. Then, there exists a positive constant $D>0$ such that
     $$\int_{\Omega}|u|^{q(x)}dx\leqslant D\left( \displaystyle\int_{Q}\varPhi_{j_1}\left(|D^su| \right)d\mu+\displaystyle\int_{Q}\varPhi_{j_2}\left(|D^su| \right)d\mu\right) 
     $$
     for all $u\in C^\infty_0(\Omega).$
       \end{lem}  
\begin{proof}[\textbf{Proof}]
First, we point out that for any $x\in \Omega$, the following inequality hold true
$$|u(x)|^{q(x)}\leqslant |u(x)|^{q^-}+|u(x)|^{q^+}~~\text{for all}~~ u\in C^\infty_0(\Omega).$$
Integrating the above inequality with request to $x\in \Omega$, we get  
$$\int_\Omega|u(x)|^{q(x)}dx\leqslant \int_\Omega|u(x)|^{q^+}dx+\int_\Omega|u(x)|^{q^-}dx~~\text{for all}~~ u\in C^\infty_0(\Omega).$$
By Poincar\'{e} inequality, there exists $c_1>0$ such that
$$
\begin{aligned}
\int_\Omega|u(x)|^{q(x)}dx &\leqslant c_1\left( \int_Q |D^su|^{q^+}d\mu+\int_Q|D^s u|^{q^-}d\mu\right)\\
& = c_1\left( \int_Q |D^su|^{\varphi^+_{j_1}}d\mu+\int_Q|D^s u|^{\varphi^-_{j_2}}d\mu\right).
\end{aligned}
$$
On the other hand, by a variant of \cite[Lemma 3]{ra2}, we infer that there exists a positive constant $c_2>0$ such that
$$\int_Q |D^su|^{\varphi^+_{j_1}}d\mu+\int_Q|D^s u|^{\varphi^-_{j_2}}d\mu\leqslant c_2\left( \int_Q \varPhi_{j_1}(|D^su|)d\mu+\int_Q\varPhi_{j_2}(|D^s u|)d\mu\right). $$
Combining the last two inequalities, we obtain the conclusion of the lemma. 
\end{proof}
\begin{lem}\label{anlem7}
Let $\lambda>0$ be fixed. Assume that the hypothesis of Theorem $\ref{anth4}$ are fulfilled. Then the following relation holds true
$$\lim\limits_{\|u\|_{\overrightarrow{\varPhi}}\rightarrow \infty}T_\lambda(u)=\infty.$$ 
\end{lem}
\begin{proof}[\textbf{Proof}]
First, we show that
\begin{equation}\label{an22}
\lim\limits_{||u||_{\overrightarrow{\varPhi}}\rightarrow \infty}\dfrac{J(u)}{I(u)}=\infty.
\end{equation}
Assume by contradiction that the above relation does not hold true. Then there exists an $M>0$ such that for each $n\in \mathbb{N}^*$, there exists a $u_n\in W^{s}_0L_{\overrightarrow{\varPhi}}(Q)$ with $||u||_{\overrightarrow{\varPhi}}>n$ and 
 \begin{equation}\label{an23}
 \dfrac{J(u_n)}{I(u_n)}\leqslant M.
 \end{equation}
 While $||u_n||_{\overrightarrow{\varPhi}}=\sum_{i=1}^{N}[u_n]_{i}\rightarrow \infty$ as $n\rightarrow \infty,$ the sequence $\left\lbrace [u_n]_{k}\right\rbrace_n$, where $k$ is given in $(\ref{an7})$ either bounded or unbounded.
 On the other hand, it is not difficult to see that 
\begin{equation}\label{+}
 \begin{aligned}
 \int_\Omega|u(x)|^{q(x)}dx & \leqslant \int_\Omega|u(x)|^{q^+}dx+\int_\Omega|u(x)|^{q^-}dx\\
 &=\|u\|^{q^+}_{q^+}+\|u\|^{q^-}_{q^-}\\
 &\leqslant c_1\left( \|u\|^{q^+}_{\varPhi_k}+\|u\|^{q^-}_{\varPhi_k}\right)\\
 &\leqslant   c_1c_2\left( [u]^{q^+}_{k}+[u]^{q^-}_{k}\right)
  \end{aligned}
\end{equation}
 for all $u\in W^{s}_0L_{\overrightarrow{\varPhi}}(Q)$. If $\left\lbrace [u_n]_{k}\right\rbrace_n$ is bounded, then by $(\ref{+})$, we have that
$\left\lbrace I(u_n)\right\rbrace_n$ is bounded. On the other hand, by $(\ref{-})$, we have 
$$J(u_n) \geqslant \dfrac{\|u_n\|^{\varphi^-_{\min}}_{\overrightarrow{\varPhi}}}{N^{\varphi_{\min}^--1}}-N.$$
Consequently, in this case we obtain that 
$$\lim\limits_{n\rightarrow \infty}\dfrac{J(u_n)}{I(u_n)}=\infty$$
which contradicts with $(\ref{an23})$.\\
Now, we assume that $ [u_n]_{k}\rightarrow \infty$ as $n\rightarrow \infty$, we can assume that $ [u_n]_{k}>1$ for all $n$. So,
$$\dfrac{J(u_n)}{I(u_n)}\geqslant \dfrac{[u_n]_{k}^{\varphi^-_k}}{c_1c_2\left( [u_n]^{q^+}_{k}+[u_n]^{q^-}_{k}\right)}.$$
Since $\varphi^-_k>q^+$, the above inequality implies that 
$$\lim\limits_{n\rightarrow \infty}\dfrac{J(u_n)}{I(u_n)}=\infty$$
which contradicts with $(\ref{an23})$.\\
Next, Assume by contradiction that the conclusion of Lemma $\ref{anlem7}$
is not valid. Then there exists an $M_1>0$, such that for each $n\in \N^*$, there exists $v_n\in W^{s}_0L_{\overrightarrow{\varPhi}}(Q)$, $\|v_n\|_{\overrightarrow{\varPhi}}>n$ and 
$$|T_\lambda(v_n)|=|J(v_n)-\lambda I(v_n)|\leqslant M_1.$$
Thus, it is clear that $\|v_n\|_{\overrightarrow{\varPhi}}\rightarrow \infty$ as $n\rightarrow \infty$ and since 
$$J(v_n) \geqslant \dfrac{\|v_n\|^{\varphi^-_{\min}}_{\overrightarrow{\varPhi}}}{N^{\varphi_{\min}^--1}}-N$$
it follows that $J(v_n)\rightarrow \infty$
 as $n\rightarrow \infty$. Then, we find that for each $n$ large enough, we have 
 $$\left| 1-\dfrac{J(v_n)}{I(v_n)}\right| \leqslant \dfrac{M_1}{J(v_n)}.$$
 Then, passing to the limit as $n\rightarrow \infty$ in the above inequality and taking into account the facts that $\dfrac{J(v_n)}{I(v_n)}\rightarrow \infty$ as $n\rightarrow \infty$ and $J(v_n)\rightarrow \infty$
  as $n\rightarrow \infty$, we obtain a contradiction.
\end{proof}                                           \subsection{Proof of mains results}
  \begin{proof}[\textbf{Proof of Theorem $\ref{anth1}$}]
  By Lemmas $\ref{anlem1}$ and $\ref{anlem2}$, we can using the mountain pass theorem $\ref{an2.2}$, and we deduce the existence of a sequence $\left\lbrace u_n\right\rbrace \in W^{s}_0L_{\overrightarrow{\varPhi}}(Q)$ such that
  \begin{equation}\label{an29}
  T_\lambda(u_n)\rightarrow c_1~~\text{and}~~ T'_\lambda(u_n)\rightarrow 0~~\text{in}~~(W^{s}_0L_{\overrightarrow{\varPhi}}(Q))^*~~\text{as}~~n\rightarrow \infty.
  \end{equation}
 We assume that $\left\lbrace u_n\right\rbrace$ is bounded in  $W^{s}_0L_{\overrightarrow{\varPhi}}(Q)$. By contradiction, we suppose that there exists a subsequence still denoted by $\left\lbrace u_n\right\rbrace $ such that $\|u_n\|_{\overrightarrow{\varPhi}}\rightarrow \infty$ and that $\|u_n\|_{\overrightarrow{\varPhi}}>1$ for all $n$. Relation $(\ref{an29})$ and the above considerations implies that for $n$ large enough, we have
 $$
 \begin{aligned}
 1+c_1+\|u_n\|_{\overrightarrow{\varPhi}}& \geqslant T_\lambda(u_n)-\dfrac{1}{q^-}\left\langle T'_\lambda(u_n), u_n\right\rangle\\
 &\geqslant \sum_{i=1}^{N} \int_{Q}\left( \varPhi_i\left(|D^su_n| \right)-\dfrac{1}{q^-}\varphi_i(|D^su_n|)D^su_n \right)d\mu  \\
 &\geqslant \left( 1-\dfrac{\varphi^+_{\max}}{q-}\right) \sum_{i=1}^{N} \int_{Q} \varPhi_i\left(|D^su_n| \right)d\mu\\
 &\geqslant
 \left( 1-\dfrac{\varphi^+_{\max}}{q-}\right) \left( \dfrac{\|u_n\|^{\varphi^-_{\min}}_{\overrightarrow{\varPhi}}}{N^{\varphi_{\min}^--1}}-N\right).
  \end{aligned}
 $$
 Since $q^->\varphi^+_{\max}$, and dividing by $\|u_n\|^{\varphi^-_{\min}}_{\overrightarrow{\varPhi}}$ in the above inequality and passing to the limit as $n\rightarrow \infty$, we obtain a contradiction. Then $\left\lbrace u_n\right\rbrace $ is bounded in $W^{s}_0L_{\overrightarrow{\varPhi}}(Q)$. This information combined with the fact that $W^{s}_0L_{\overrightarrow{\varPhi}}(Q)$ is reflexive, implies that there exists a subsequence still denoted by $\left\lbrace u_n\right\rbrace $ and $u_0\in W^{s}_0L_{\overrightarrow{\varPhi}}(Q)$ such that $\left\lbrace u_n\right\rbrace $ converges weakly to $u_0$ in $W^{s}_0L_{\overrightarrow{\varPhi}}(Q)$. On the other hand, since $W^{s}_0L_{\overrightarrow{\varPhi}}(Q)$ is compactly embedded in $L^{q(x)}(\Omega)$, it follows that $\left\lbrace u_n\right\rbrace $ converges strongly to $u_0$ in $L^{q(x)}(\Omega)$.
 Then by H\"{o}lder inequality, we have that 
 $$ \lim\limits_{n\rightarrow \infty}\int_\Omega |u_n|^{q(x)-2}u_n(u_n-u_0)dx=0.$$
 This fact and relation $(\ref{an29})$, implies that
 $$\lim\limits_{n\rightarrow \infty}\left\langle T'_\lambda(u_n), u_n-u_0\right\rangle =0.$$
 Thus we deduce that 
 \begin{equation}\label{an31}
 \lim\limits_{n\rightarrow \infty}\sum_{i=1}^{N} \int_{Q}a_i(|D^su_n|)D^su_n\left( D^su_n-D^su_0 \right)d\mu =0.
 \end{equation}
 Since $\left\lbrace u_n\right\rbrace $ converge weakly to $u_0$ in $W^{s}_0L_{\overrightarrow{\varPhi}}(Q)$, by relation $(\ref{an31})$, we find that 
 {\small \begin{equation}\label{an32}
  \lim\limits_{n\rightarrow \infty}\sum_{i=1}^{N} \int_{Q}\left( a_i(|D^su_n|)D^su_n-a_i(|D^su_0|)D^su_0\right) \left( D^su_n-D^su_0 \right)d\mu =0.
 \end{equation}}
 Since, for each $i\in \left\lbrace 1,...,N\right\rbrace $, $\varPhi_i$ is convex, we have
 $$\varPhi_i(|D^su|)\leqslant \varPhi_i\left( \dfrac{|D^su+D^sv|}{2}\right) +a_i(|D^su|)D^su\dfrac{D^su-D^sv}{2}$$
  $$\varPhi_i(|D^sv|)\leqslant \varPhi_i\left( \dfrac{|D^su+D^sv|}{2}\right) +a_i(|D^sv|)D^sv\dfrac{D^sv-D^su}{2}$$
  for every $u,v \in W^{s}_0L_{\overrightarrow{\varPhi}}(Q)$.
  Adding the above two relations and integrating over $Q$, we find that
\begin{equation}\label{an33}
  \begin{aligned}
  \dfrac{1}{2}&\int_{Q} \left( a_i(|D^su|)D^su-a_i(|D^sv|)D^sv\right) \left( D^su-D^sv \right)d\mu\\
  &\geqslant \int_Q\varPhi_i(|D^su|) d\mu+ \int_Q\varPhi_i(|D^sv|) d\mu- 2\int_Q\varPhi_i\left( \dfrac{|D^su-D^sv|}{2}\right)d\mu,
    \end{aligned}
\end{equation}
 for every $u,v \in W^{s}_0L_{\overrightarrow{\varPhi}}(Q)$, and each $i\in \left\lbrace 1,...,N\right\rbrace.$
 
 On the other hand, since for each $i\in \left\lbrace 1,...,N\right\rbrace $ we know that $\varPhi_i~:~ [0,\infty) \rightarrow \R$ is an increasing continuous function, with $\varPhi_i(0)=0$. Then by the conditions $(\ref{f1.})$ and $(\ref{f2.})$, we can apply \cite[Lemma 2.1]{Lam}  in order to obtain
 \begin{equation}\label{an34}
 \begin{aligned}
 \dfrac{1}{2}&\left[ \int_Q\varPhi_i(|D^su|) d\mu+ \int_Q\varPhi_i(|D^sv|) d\mu\right] \\
& \geqslant \int_Q\varPhi_i\left( \dfrac{|D^su+D^sv|}{2}\right)d\mu+\int_Q\varPhi_i\left( \dfrac{|D^su-D^sv|}{2}\right)d\mu,
 \end{aligned}
 \end{equation}
 for every $u,v \in W^{s}_0L_{\overrightarrow{\varPhi}}(Q)$, and each $i\in \left\lbrace 1,...,N\right\rbrace.$ 
 By $(\ref{an33})$ and $(\ref{an34})$, it follows that for each $i\in \left\lbrace 1,...,N\right\rbrace$, we have
 \begin{equation}\label{an35}
 \begin{aligned}
 & \int_{Q} \left( a_i(|D^su|)D^su-a_i(|D^sv|)D^sv\right) \left( D^su-D^sv \right)d\mu\\
 &\geqslant 4\int_Q\varPhi_i\left( \dfrac{|D^su-D^sv|}{2}\right)d\mu
 \end{aligned}
 \end{equation}
 for every $u,v \in W^{s}_0L_{\overrightarrow{\varPhi}}(Q)$, and each $i\in \left\lbrace 1,...,N\right\rbrace.$ \\
 Relations $(\ref{an32})$ and $(\ref{an35})$ show that $\left\lbrace u_n\right\rbrace $ converge strongly to $u_0$ in $W^{s}_0L_{\overrightarrow{\varPhi}}(Q)$. Then by relation $(\ref{an29})$, we have 
 $T_\lambda(u_0)=c_1>0$ and $T'_\lambda(u_0)=0$, that is, $u_0$ is a non trivial weak solution.
    \end{proof}
  \begin{proof}[\textbf{Proof of Theorem $\ref{anth2}$}]    
   Let $\lambda_*>0$ be defined as in $(\ref{an40})$ and $\lambda\in (0,\lambda_*)$. By Lemma $\ref{anlem3}$ it follows that on the boundary oh the ball centered in the origin and of radius $\rho$ in $W^{s}_0L_{\overrightarrow{\varPhi}}(Q)$, denoted by $B_\rho(0)$, we have 
         $$\inf\limits_{\partial B_\rho(0)}T_\lambda>0.$$
      On the other hand, by Lemma $\ref{anlem4}$, there exists $\phi \in W^{s}_0L_{\overrightarrow{\varPhi}}(Q)$ such that $T_\lambda(t\phi)<0$ for all $t>0$ small enough. Moreover for any $u\in B_\rho(0)$, we have 
      $$
               \begin{aligned}
           T_\lambda(u)\geqslant \dfrac{1}{N^{\varphi_{\min}^--1}} \|u\|^{\varphi^-_{\min}}_{\overrightarrow{\varPhi}}-\dfrac{\lambda c_1^{q^-}}{q^-}\|u\|^{q^-}_{\overrightarrow{\varPhi}}.
                 \end{aligned}
                            $$
      It follows that
      $$-\infty<c:=\inf\limits_{\overline{B_\rho(0)}} T_\lambda<0.$$   
      We let now $0<\varepsilon <\inf\limits_{\partial  B_\rho(0)}  T_\lambda -  \inf\limits_{B_\rho(0)} T_\lambda.$    Applying Theorem $\ref{anek}$ to the functional 
      $T_\lambda : \overline{B_\rho(0)}\longrightarrow \R$, we find $u_\varepsilon \in \overline{B_\rho(0)}$ such that 
       $$
            \left\{ 
                 \begin{array}{clclc}
               T_\lambda(u_\varepsilon)&<\inf\limits_{\overline{B_\rho(0)}} T_\lambda+\varepsilon,& \\\\
                 T_\lambda(u_\varepsilon)&< T_\lambda(u)+\varepsilon ||u-u_\varepsilon||_{\overrightarrow{\varPhi}},& \text{  } u\neq u_\varepsilon.
                 \end{array}
                 \right. 
              $$
       Since  $T_\lambda(u_\varepsilon)\leqslant  \inf\limits_{\overline{B_\rho(0)}} T_\lambda+\varepsilon\leqslant \inf\limits_{B_\rho(0)} T_\lambda+\varepsilon < \inf\limits_{\partial  B_\rho(0)}  T_\lambda$, we deduce $u_\varepsilon  \in B_\rho(0)$. 
       
       Now, we define $\Lambda_\lambda :  \overline{B_\rho(0)}\longrightarrow \R$ by 
       $$\Lambda_\lambda(u)=T_\lambda(u)+\varepsilon||u-u_\varepsilon||_{\overrightarrow{\varPhi}}.$$
       It's clear that $u_\varepsilon$ is a minimum point of $\Lambda_\lambda$ and then
       $$\dfrac{\Lambda_\lambda(u_\varepsilon+t v)-\Lambda_\lambda(u_\varepsilon)}{t}\geqslant 0$$ 
       for small $t>0$, and any $v\in B_\rho(0).$ The above relation yields 
           $$\dfrac{T_\lambda(u_\varepsilon+t v)-T_\lambda(u_\varepsilon)}{t}+\varepsilon||v||_{\overrightarrow{\varPhi}}\geqslant 0.$$
           Letting $t\rightarrow$ it follows that $\left\langle T'_{\lambda}(u_\varepsilon),v\right\rangle +\varepsilon ||v||_{\overrightarrow{\varPhi}}>0$ and we infer that $$||T'_{\lambda}(u_\varepsilon)||_{\overrightarrow{\varPhi},*}\leqslant \varepsilon.$$
            We deduce that there exists a sequence $\left\lbrace v_n\right\rbrace \subset B_\rho(0)$ such that 
            \begin{equation}\label{an10}
            T_\lambda(v_n) \longrightarrow c \text{ and } T'_\lambda(v_n)\longrightarrow 0.
            \end{equation}
         It is clear that $\left\lbrace v_n\right\rbrace $ is bounded in $W^{s}_0L_{\overrightarrow{\varPhi}}(Q)$. Thus, there exists $v\in W^{s}_0L_{\overrightarrow{\varPhi}}(Q)$, such that up to a subsequence   $\left\lbrace v_n\right\rbrace $ converges weakly to $v$ in $W^{s}_0L_{\overrightarrow{\varPhi}}(Q)$. Actually with similar arguments to those used of the end of Theorem $\ref{anth1}$, we can show that $\left\lbrace v_n\right\rbrace $ is converges strongly to $v$ in $W^{s}_0L_{\overrightarrow{\varPhi}}(Q)$. Thus, by   $(\ref{an10})$
         $$T_\lambda(v)=c <0~~\text{ and }~~ T'_\lambda(v)=0.$$
         Then, $v$ is a nontrivial weak solution for Problem \hyperref[P]{$(P_a)$}. This complete the proof.
             \end{proof}
     
       \begin{proof}[\textbf{Proof of Theorem $\ref{anth3}$}]   
    The existence of a positive constant $\lambda_*$ such that any $\lambda\in (0,\lambda_*)$ is an eigenvalue of Problem \hyperref[P]{$(P_a)$} is an immediate  consequence of Theorem $\ref{anth2}$. In order to prove the second part of Theorem $\ref{anth3}$, we will show for $\lambda$ positive and large enough the functional $T_\lambda$ possesses a nontrivial global minimum point in  $W^{s}_0L_{\overrightarrow{\varPhi}}(Q)$. For that,  by Theorem $\ref{andr}$  it is enough  to prove that $T_\lambda$ is coercive and   weakly lower semi continuous. 
  
 First, we show that $T_\lambda$ is weakly lower semi continuous. Let $\left\lbrace u_n\right\rbrace \subset W^{s}_0L_{\overrightarrow{\varPhi}}(Q)$ be a sequence which convergences weakly to $u$ in $ W^{s}_0L_{\overrightarrow{\varPhi}}(Q)$. Since $J$ is weakly lower semi continuous, then
     $$J(u)\leqslant \liminf\limits_{n\rightarrow \infty} J(u_n).$$
     On the other hand, since $W^{s}_0L_{\overrightarrow{\varPhi}}(Q)$ is compactly embedded in $L^{q(x)}(\Omega)$. Then, $u_n$ is converges
     strongly to $u$ in $L^{q(x)}(\Omega)$  and by Proposition $\ref{anproop5}$, we have
     $$\lim\limits_{n\rightarrow \infty}I(u_n)=I(u).$$
    Thus, we find 
     $$T_\lambda(u)\leqslant \liminf\limits_{n\rightarrow \infty} T_\lambda(u_n).$$
     Therefore, $T_\lambda$ is weakly lower semi continuous. Next, by Lemma $\ref{anlem5}$, the functional $T_\lambda$ is also coercive on   $W^{s}_0L_{\overrightarrow{\varPhi}}(Q)$. These two facts enable us to apply Theorem $\ref{andr}$ in order to find that there exists $u_\lambda \in  W^{s}_0L_{\overrightarrow{\varPhi}}(Q)$. a global minimizer of $T_\lambda$, and thus a weak solution of Problem \hyperref[P]{$(P_a)$}.
 
Now, we show that $u_\lambda$ is non trivial. Indeed, letting $t_0>1$ be fixed real and $u_0(x)=t_0$ for all $x\in \Omega$, we have $ u_0\in W^{s}_0L_{\overrightarrow{\varPhi}}(Q)$ and 
$$
\begin{aligned}
T_\lambda(u_0)&=\displaystyle\int_{Q}\sum_{i=1}^{N}\varPhi_i\left(|D^su_0| \right)d\mu-\lambda\displaystyle\int_{\Omega}\dfrac{1}{q(x)}|u_0|^{q(x)}dx\\
&\leqslant L-\dfrac{\lambda}{q^+}\int_{\Omega} |u_0|^{q(x)} dx\\
&\leqslant L-\dfrac{\lambda}{q^+}|t_0|^{q^-}|\Omega|
\end{aligned}
$$
where $L$ is a positive constant. Thus, for $\lambda^*>0$ large enough, $T_\lambda(u_0)<0$ for any $\lambda\in [\lambda^*,\infty)$. It follows that $T_\lambda(u_\lambda)<0$ for any $\lambda\in [\lambda^*,\infty)$ and thus $u_\lambda$ is a nontrivial weak solution of Problem  \hyperref[P]{$(P_a)$} for any $\lambda\in [\lambda^*,\infty)$. Therefore, Problem  \hyperref[P]{$(P_a)$} has a nontrivial weak solution for all $\lambda\in(0,\lambda_*)\cup [\lambda^*,\infty)$.
                  \end{proof}    
       \begin{proof}[\textbf{Proof of Theorem $\ref{anth4}$}] 
   \text{  }  \\
      
   $ \bullet  $ \textbf{Step 1.} We note that by Lemma $\ref{anlem6}$, we can conclude that
   $$\lambda_1>0.$$
   On the other hand,  We point out by definition of $\varphi_i^-$, for any $i=1,...,N$ that 
   $$a_i(t).t^2=\varphi_i(t)t\geqslant \varphi_i^-\varPhi_i(t)~~\forall t\geqslant 0~~i=1,...,N.$$
   The above inequality and Lemma $\ref{anlem6}$ imply 
   \begin{equation}\label{e21}
   \lambda_0:=\inf\limits_{u\in W^{s}_0L_{\overrightarrow{\varPhi}}(Q)\setminus\left\lbrace 0\right\rbrace }\dfrac{J_1(u)}{I_1(u)}>0.\end{equation}

\textbf{$\bullet$ Step 2.}  We show that any $\lambda\in (0,\lambda_0)$ is not an eigenvalue of Problem \hyperref[P]{$(P_a)$}, where $\lambda_0$ is given by $(\ref{e21})$. Indeed, assuming by contradiction that there exists $\lambda\in (0,\lambda_0)$ an eigenvalue of Problem \hyperref[P]{$(P_a)$}. That is, it follows that there exists $u_\lambda\in W^{s}_0L_{\overrightarrow{\varPhi}}(Q)\setminus\left\lbrace 0\right\rbrace$ such that 
 $$\left\langle J'(u_\lambda),v\right\rangle =\lambda \left\langle I'(u_\lambda),v\right\rangle ~~\forall v\in W^{s}_0L_{\overrightarrow{\varPhi}}(Q).$$
 Then, for $v=u_\lambda$ we have
 $$\left\langle J'(u_\lambda),u_\lambda\right\rangle =\lambda \left\langle I'(u_\lambda),u_\lambda\right\rangle,$$
 so,  
 $$J_1(u_\lambda)=\lambda I_1(u_\lambda).$$
 The fact that $u_\lambda\in W^{s}_0L_{\overrightarrow{\varPhi}}(Q)\setminus\left\lbrace 0\right\rbrace$ implies that $I_1(u_\lambda)>0$. Since $\lambda<\lambda_0$, then 
 $$J_1(u_\lambda)\geqslant \lambda_0 I_1(u_\lambda)>\lambda I_1(u_\lambda)=J_1(u_\lambda).$$
 Clearly, the above inequality lead, to a contradiction. Thus, step $2$ is verified.
 
\textbf{$\bullet$ Step 3} We show that any $\lambda\in (\lambda_1,\infty)$ is an eigenvalue of Problem \hyperref[P]{$(P_a)$}.\\
  Let $\lambda\in (\lambda_1,\infty)$ be arbitrary.  Thus, $\lambda$ is an eigenvalue of Problem \hyperref[P]{$(P_a)$} if and only if there exists $u_\lambda\in W^{s}_0L_{\overrightarrow{\varPhi}}(Q)\setminus\left\lbrace 0\right\rbrace $ a critical point of $T_\lambda$.  
 Indeed, By Lemma $\ref{anlem7}$, we can obtain that $T_\lambda$ is coercive, and as  $T_\lambda$ is weakly lower semi continuous. Then by apply Theorem $\ref{andr}$ in order prove that there exists $u_\lambda \in  W^{s}_0L_{\overrightarrow{\varPhi}}(Q)$ a global minimizing point of $T_\lambda$, and thus a critical point of $T_\lambda$. Next we prove that $u_\lambda\neq 0$. Indeed since
    $\lambda_1=\inf\limits_{w\in W^{s}_0L_{\overrightarrow{\varPhi}}(Q)\setminus\left\lbrace 0\right\rbrace }\dfrac{J(w)}{I(w)}$ and $\lambda>\lambda_1$, then there exists $v_\lambda \in W^{s}_0L_{\overrightarrow{\varPhi}}(Q)$ such that
    $$J(v_\lambda)<\lambda I(v_\lambda)$$
    this implies that 
    $$T_\lambda(v_\lambda)<0.$$
    Thus
    $$\inf\limits_{w\in W^{s}_0L_{\overrightarrow{\varPhi}}(Q)} T_\lambda<0.$$
    So, we conclude that $u_\lambda$ is a nontrivial critical point of $T_\lambda$, that is, $\lambda$ is an eigenvalue of Problem \hyperref[P]{$(P_a)$}. Then, step 3 is verified.

 \textbf{$\bullet$}
Finally by the above steps, we deduce that $\lambda_0\leqslant \lambda_1$. The proof of Theorem $\ref{anth4}$ is now completed.
                     \end{proof}  
                     
  \subsection*{Disclosure statement}
             No potential conflict of interest was reported by the author.

\end{document}